%
%
%

\input amstex
\documentstyle{amsppt}
\magnification=1200
 \vsize19.5cm
  \hsize13.5cm \TagsOnRight
\pageno=1
\baselineskip=15.0pt
\parskip=3pt

\def\p{\partial}
\def\noo{\noindent}
\def\eps{\varepsilon}

\def\Om{\Omega}

\def\pom{{\p \Om}}
\def\bom{{\overline\Om}}
\def\R{\bold R}

\def\wtt{\tilde}

\def\Ga{\Gamma}

\def\ol{\overline}
\def\ul{\underline}

\def\phi{\varphi}

\nologo
 \NoRunningHeads

\topmatter

\title{The mean curvature measure}\endtitle

\author{Qiuyi Dai\ \ \  Neil S. Trudinger\ \ \ Xu-Jia Wang }\endauthor

\address
\newline
Qiuyi Dai: Department of Mathematics, Hunan Normal University,
Changsha, 410081, P.R. China.
\newline
Neil S. Trudinger and Xu-Jia Wang:
 Centre for Mathematics and its Applications, Australian National
    University,  Canberra ACT 0200, Australia
 \endaddress

\email \newline
 daiqiuyi\@yahoo.com.cn, \ \
 neil.trudinger\@maths.anu.edu.au, \ \
 wang\@maths.anu.edu.au.
\endemail

\thanks{
This work was supported by ARC grants DP0664517 and DP0879422; and
NSFC grants 10428103 and 10671064. }\endthanks

\abstract {We assign a measure to an upper semicontinuous function
which is subharmonic with respect to the mean curvature operator,
so that it agrees with the mean curvature of its graph when the
function is smooth. We prove that the measure is weakly continuous
with respect to almost everywhere convergence. We also establish a
sharp Harnack inequality for the minimal surface equation, which
is crucial for our proof of the weak continuity. As an application
we prove the existence of weak solutions to the corresponding
Dirichlet problem when the inhomogeneous term is a measure.
}\endabstract

\endtopmatter


\document

\baselineskip=13.8pt
\parskip=3pt

\centerline {\bf 1. Introduction}

\vskip10pt

Notions of curvature measures arise in convex geometry, (see for
example [S]), and were extended to general surfaces by Federer
[F1] under a hypothesis of {\it positive reach}. For graphs of
functions, this condition is equivalent to semi-convexity and
implies twice almost everywhere differentiability by virtue of the
well-know theorem of Aleksandrov. The development of a
corresponding theory of curvature measures on more general sets is
an open problem. Without any assumption such a theory seems
impossible as the second derivative of a nonsmooth function is
usually a distribution but not a measure. In this paper we
consider the {\it mean curvature} and restrict ourselves to graphs
of functions defined over domains $\Om$ in Euclidean $n$-space,
$\R^n$. The mean curvature has been the most extensively studied
geometric quantity but usually it is regarded as a distribution
when the function is not twice differentiable, such as in the case
when its graph is a rectifiable set.

In particular in this paper we assign a measure to an upper
semicontinuous function which is subharmonic with respect to the
mean curvature operator, so that it agrees with the mean curvature
of its graph when the function is smooth. We prove that the
measure is weakly continuous with respect to almost everywhere
convergence (Theorem 6.1). We also establish a sharp Harnack
inequality for the minimal surface equation (Theorem 2.1), which
is crucial for our proof of weak continuity. As an application we
prove the existence of weak solutions to the Dirichlet problem of
the mean curvature equation when the right hand side is a measure
(Theorem 7.1).

We say an upper semi-continuous function $u: \Om \rightarrow
[-\infty, +\infty)$ is subharmonic with respect to the mean
curvature operator $H_1$, or {\it $H_1$-subharmonic} in short, if
the set $\{u=-\infty\}$ has measure zero and $H_1[u]\ge 0$ in the
viscosity sense. That is for any open set $\omega\subset\Om$ and
any smooth function $h\in C^2(\ol \omega)$ with $H_1[h]\le 0$,
$h\geq u$ on $\p\omega$, one has $h\geq u$ in $\omega$. We say a
function $u$ is {\it $H_1$-harmonic} if it is $H_1$-subharmonic
and for any open set $\omega\subset\Om$ and any $H_1$-subharmonic
function $h$ in $\omega$ with $h\le u$ on $\p\omega$, one has
$h\le u$ in $\omega$. This definition does not imply directly that
an $H_1$-harmonic function is bounded from below, but we will
prove in Section 4 it is the case, and so is smooth.  We denote
the set of all $H_1$-subharmonic functions in $\Om$ by
$SH_1(\Om)$.

A main result of the paper is the weak continuity of the mean
curvature operator. That is if $\{u_k\}$ is a sequence of smooth
$H_1$-subharmonic functions which converges a.e. to $u\in
SH_1(\Om)$, then $H_1[u_k]$ converges weakly to the density of a
measure $\mu$. The measure $\mu$ depends only on $u$ but not on
the sequence $\{u_k\}$, so that we can assign a measure, called
{\it the mean curvature measure} and denoted by $\mu_1[u]$, to the
function $u$. Note that our measure $\mu_1$ is defined on $\Om$
but Federer's measure $\nu_1$ is defined on the graph of $u$.

A crucial ingredient for the proof of the weak continuity is a
refined Harnack inequality, also established in this paper,  for
the minimal surface equation
$$H_1[u]=: \text{div}(\frac{Du}{\sqrt{1+|Du|^2}})=0.\tag 1.1$$
Namely
$$\sup_{B_r} u\le C\inf_{B_r} u\tag 1.2$$
for nonnegative solution of (1.1) in $\ol B_{2r}$. The Harnack
inequality for the mean curvature equation has been studied in
several works [FL, Lia, PS1, T1]. We prove that the  constant $C$
depends on the decay rate of $|\{x\in B_{2r}:\ u(x)>t\}|_n$, or
$|\{x\in \p B_{2r}:\ u(x)>t\}|_{n-1}$, as $t\to\infty$, where
$|\cdot |_k$ denotes the $k$-dimensional Haudorff measure. This is
indeed the best possibility one can expect. A similar Harnack
inequality also holds for the non-homogeneous equation, see Remark
2.4.

As an application, we study the existence of solutions to the
Dirichlet problem of the mean curvature equation
$$\align
  H_1[u] & =\nu\ \ \ \ \text{in}\ \ \Om, \tag 1.3\\
  u  & =  \phi \ \ \ \  \text{on}\ \ \pom,\\
  \endalign $$
where $\nu$ is the density of a nonnegative measure, with respect
to Lebesgue measure.

For the Dirichlet problem of the mean curvature equation, it is
usually assumed that the right hand side $\nu$ is a Lipschitz
function, so that the interior gradient estimate holds and the
solution is smooth, in $C^{2, \alpha}(\Om)$ for $\alpha\in (0, 1)$
[GT]. If $\nu$ is not Lipschitz continuous, the solution may not be
$C^2$ smooth even if $\nu$ is H\"older continuous; (see the example
in \S 8). In [Gia, G2] it was proved that when $\nu$ is a measurable
function satisfying a necessary condition, equation (1.3) has a weak
solution which is a minimizer of an associated functional.  Through
the mean curvature measure introduced above, we introduce a notion
of weak solution and prove its existence when $\nu$ is a nonnegative
measure.

This paper is arranged as follows. In Section 2 we establish the
Harnack inequality (1.2) for the minimal surface equation. In
Section 3 we establish an integral gradient estimate and a uniform
estimate for $H_1$-subharmonic functions. In these two sections we
assume that the functions are smooth. But the assumption can be
removed by an approximation result proved in Section 5.

In Section 4 we introduce the Perron lifting and prove some basic
properties for $H_1$-subharmonic functions. In Section 5 we prove
that every $H_1$-subharmonic function can be approximated by a
sequence of smooth, $H_1$-subharmonic functions. Section 6 is
devoted to the proof of the weak continuity of the mean curvature
operator. The Dirichlet problem is discussed in Section 7. The
final Section 8 contains some remarks.

In recent years it was proved that for several important
homogeneous elliptic operators, such as the $p$-Laplace operator
and the $k$-Hessian operator, one can assign a measure to a
function which is subharmonic with respect to the operators, and
as applications various potential theoretical results have been
established. See [HKM, Lab, TW1-TW4]. Our treatment of the weak
continuity of the mean curvature operator was inspired by the
earlier works [TW1-TW4]. However as the mean curvature operator is
non-homogeneous, the situation is much more delicate.

\vskip20pt

\centerline{\bf 2. The Harnack inequality}

\vskip10pt

In this section we prove a Harnack inequality for the minimal
surface equation, which will be used for the Perron liftings
process in Section 4 and the study of the Dirichlet problem in
Section 7. We also establish a weak Harnack inequality for
$H_1$-subharmonic functions, which will be used in the proof of
Lemma 4.3.

First we quote the basic existence and regularity result for the
mean curvature equation [GT]. The regularity of the mean curvature
equation is based on the interior gradient estimate (see Theorem
16.5 in [GT]).

\proclaim{Lemma 2.1} Let $u\le 0$ be a $C^3$ solution to the mean
curvature equation
$$H_1[u]=f(x)\ \ \text{in}\ B_r(0).\tag 2.1$$
Then
$$ |Du(0)|\le C_1e^{C_2\frac{|u(0)|}{r}},\tag 2.2$$
where $C_1, C_2$ depend only on $n$ and $\|f\|_{C^{0,1}}$.
\endproclaim

Simpler proofs of the interior gradient estimate, with
$\frac{|u(0)|}{r}$ replaced by $\frac{|u(0)|^2}{r^2}$, was given
in [K1, Wan]. The proofs also applies to the $k$-th mean curvature
equation and more general Weingarten curvature equations [K2,
Wan].

From the gradient estimate, the mean curvature equation becomes
uniformly elliptic and one has local uniform estimate in $C^{2,
\alpha}$ for the equation, for any $\alpha\in (0, 1)$.

By the regularity, one has the existence of solutions to the
Dirichlet problem (see Theorem 16.8 in [GT]).

\proclaim{Lemma 2.2}
 Let $\Om$ be a bounded smooth domain in $\R^n$. Suppose the mean curvature
of $\pom$ is positive. Then for any continuous function $\phi$ on
$\pom$, there is a unique solution $u\in C^2(\Om)\cap C^0(\ol\Om)$
to $H_1[u]=0$ such that $u=\phi$ on $\pom$.
\endproclaim

Lemma 2.2 also holds for the inhomogeneous equation $H_1[u]=f$
with $f\in C^{0,1}$, under certain conditions on $f$ and $\pom$,
see Theorem 16.10 in [GT].

In this section we prove the following Harnack inequality. Here we
consider smooth solutions only. In Section 4 we will show that an
$H_1$-harmonic function must be smooth.

\proclaim {Theorem 2.1} Let $u\ge 0$ be a smooth solution to the
minimal surface equation
$$H_1[u]=0 \ \ \ \text{in}\ \ol B_r(0). \tag 2.3$$
Let
$$\psi(t)=|\{x\in\p B_r(0):\ u(x)>t\}|_{n-1}, $$
where $|\cdot|_{n-1}$ denotes the $(n-1)$-dim Hausdorff measure.
Suppose $\psi(t)\to 0$ as $t\to\infty$. Then there exists a
constant $C>0$ depending only on $n$, $r$, and $\psi$ such that
$$\sup_{B_{r/2}(0)} u\le C\inf_{B_{r/2}(0)} u. \tag 2.4$$
\endproclaim

\noo{\bf Remark 2.1}
\newline
(i) The Harnack inequality (2.4) was also established in [T1], but
the constant $C$ depends on $\sup u$. The main point in [T1] is a
positive lower bound of $u(0)$ for the mean curvature equation and
more general elliptic equations satisfying certain structural
conditions. The paper [T1] also includes the following weak Harnack
inequality for the upper bound for $u(0)$: if $u\in W^{2,
n}(B_r(0))$ is a subsolution, then for any $p\in (0, n]$,
$$\sup_{B_{r/2}} u\le \frac{C}{r^{(n+2)/p}}
       \big(\int_{B_r} (u^+)^{p+2}\big)^{1/p}, \tag 2.5$$
where $C$ is a constant depending only on $n$ and $p$. We also
refer the reader to [FL, Lia, PS1] for discussions of the Harnack
inequality.
\newline
(ii)  Recall that in the Harnack inequality for the Laplace
equation, the constant $C$ depends only on $n$. But this is
impossible for the minimal surface equation. One can construct a
positive solution of (2.3) in $B_1(0)$ such that $u(0)\le 1$ but
$\int_{B_1} u^p$ can be as large as we want, for any $p>0$. To see
this, let $\phi(x_1)$ be a positive, convex function defined for
$x_1\in (-1, 1)$ such that $\phi(x_1)$ is small when $x_1<\frac
14$ and $\phi(x_1)\to\infty$ as $x_1\to 1$. Let $u$ be the
solution of (2.3) with the Dirichlet condition $u=\phi$ on $\p
B_1$. Then by the convexity of $\phi$, $H_1[\phi]\ge 0$. Hence by
the comparison principle, we have $u\ge \phi$ in $B_1$. Hence
$\int_{B_1} u^p$ can be as large as we want provided $\phi$ is
sufficiently large near $x_1=1$. On the other hand, by
constructing a suitable upper barrier one has $u(0)\le 1$.

\vskip10pt

To prove Theorem 2.1, we start with some technical lemmas.

Let $\Om$ be an open set contained in $B_r(0)$. For $s\in (0, r]$,
denote
$$\align
\Ga^{\text{int}}_s & =\bom\cap\p B_s(0),\\
\Ga^{\text{bdy}}_s & =\pom\cap B_s(0), \\
\endalign $$
so that
$$\Ga^{\text{bdy}}_s\cup\Ga^{\text{int}}_s=\p (\Om\cap B_s(0)). $$
Let $\wtt\Ga^{\text{int}}_s$ be a geodesic ball in $\p B_s$, with
center at $(s, 0, \cdots, 0)$, such that
$|\wtt\Ga^{\text{int}}_s|_{n-1}=|\Ga^{\text{int}}_s|_{n-1}$, where
$|\cdot|_k$ denotes the $k$-dimensional Hausdorff measure. Denote by
$\rho(s)$ the geodesic radius of $\wtt\Ga^{\text{int}}_s$. Then
$$(1-\eps)\alpha_{n-1}\rho^{n-1}(s)
           \le |\Ga^{\text{int}}_s|_{n-1}\le \alpha_{n-1}\rho^{n-1}(s)$$
with $\eps\to 0$ as $\rho(s)\to 0$, where $\alpha_{n-1}$ is the
volume of the unit ball in $\R^{n-1}$. The second inequality is due
to the positive curvature of the sphere, and the first one can be
obtained easily by representing $\wtt\Ga^{\text{int}}_s$ as a graph.

Let $\rho_1$ be the constat such that
$|\Ga^{\text{int}}_s|_{n-1}\ge\frac 14 \alpha_{n-1}\rho^{n-1}(s)$
for any $\rho(s)\le \rho_1$. We also denote
$$b_n=2^{-4n}c_{n-1}\alpha_{n-1}^{-1/(n-1)}, $$
where $c_{n-1}$ is the best constant in the isoperimetric
inequality, see (2.9) below.

\proclaim{Lemma 2.3} Let $\Om$ be an open set in $B_r(0)$ for some
$\frac 14\le r\le 1$. Suppose $\Ga^{\text{bdy}}_r$ is smooth,
$\rho(r)\le \rho_1$, and
$$\rho(s)\ge \frac 14 \rho(r) \ \ \ \forall\ s\in (r', r),\tag 2.6$$
where $r'=r-\rho(r)/2b_n$. Then
$$|\Ga^{\text{bdy}}_r|_{n-1}\ge 2 |\Ga^{\text{int}}_r|_{n-1}. \tag 2.7$$
\endproclaim

\noo{\bf Proof}. We claim that
$$|\Ga^{\text{bdy}}_{r}|_{n-1}
 \ge \int_0^r |\p\Ga^{\text{int}}_{s}|_{n-2} \, ds. \tag 2.8$$
In the following we will drop the subscripts $k$ in the Hausdorff
measure $|\cdot|_k$ ($k=1, \cdots, n$) if no confusions arise.

Formula (2.8) can be derived as follows. For any point $x_0\in
\p\Ga^{\text{int}}_{r}$, by a rotation of axes we assume that
$x_0=(r, 0, \cdots, 0)$ such that $(0, \cdots, 0, 1)$ is the normal
of $\p\Ga^{\text{int}}_{r}$ at $x_0$. Then near $x_0$,
$\Ga^{\text{bdy}}_{r}$ can be represented as $x_n=\psi(x')$ such
that $\p_{x_i}\psi(x_0)=0$ for $i=2, \cdots, n-1$, where $x'=(x_1,
\cdots, x_{n-1})$. Hence at $x_0$ the area element is
$$\align
d\sigma & =\sqrt{1+|D\psi|^2}dx'\\
        & =\sqrt{1+\psi_{x_1}^2}dx' \ge dx'\\
 \endalign $$
Hence
$$\align
|\Ga^{\text{bdy}}_{r}|
 &=\int_{\Ga^{\text{bdy}}_{r}}d\sigma
 \ge \int_{\Ga^{\text{bdy}}_{r}}dx'\\
 & = \int_0^r | \p\Ga^{\text{int}}_{s}|\, ds\\
 \endalign $$
and we obtain (2.8).

By the isoperimetric inequality,
$$\align
\int_0^r | \p\Ga^{\text{int}}_{s} |\, ds
 &\ge c_{n-1} \int_0^r \big| \Ga^{\text{int}}_{s} \big|^{\frac{n-2}{n-1}}ds \tag 2.9\\
 &\ge c_{n-1}\int_{r'}^r \big|\Ga^{\text{int}}_{s} \big|^{\frac{n-2}{n-1}}ds.\\
   \endalign $$
Since $\rho(s)\ge \frac 14\rho(r)$ for any $s\in (r', r)$,
$$|\Ga^{\text{int}}_{s}|\ge 4^{-n}|\Ga^{\text{int}}_{r}|\ \ \ \forall\ s\in (r', r).$$
We obtain
$$\int_0^r | \p\Ga^{\text{int}}_{s} |\, ds
 \ge \frac {c_{n-1}}{4^n}(r-r') |\Ga^{\text{int}}_{r}|^{\frac {n-2}{n-1}}.$$
Therefore
$$\align
|\Ga^{\text{bdy}}_{r}|
  & \ge \frac {c_{n-1}}{4^n}\frac {r-r'}{|\Ga^{\text{int}}_{r}|^{1/(n-1)}} |\Ga^{\text{int}}_{r}|\\
  & \ge  \frac {c_{n-1} \alpha_{n-1}^{-1/(n-1)}}{4^{n+1}b_n} |\Ga^{\text{int}}_{r}|.\\
  \endalign $$
The Lemma holds by our choice of $b_n$. $\square$

\proclaim{Lemma 2.4} Let $\Om$ be an open set in $B_r(0)$ for some
$\frac 14\le r\le 1$. Suppose $\rho(r)\le \rho_1$ and
$\Ga^{\text{bdy}}_r$ is smooth. Then
$$|\Ga^{\text{bdy}}_r|\ge (1-4^{-n+1}-\eps) |\Ga^{\text{int}}_r|, \tag 2.10$$
with $\eps\to 0$ as $\rho(r)\to 0$. In particular there exists a
constant $\rho_2>0$ such that when $\rho(r)\le \rho_2$,
$$|\Ga^{\text{bdy}}_r|
   \ge (1- 4^{-n+5/4}) |\Ga^{\text{int}}_r| \eqno (2.10)'$$
\endproclaim

\noo{\bf Proof}. If $\rho(s)\ge \frac 14 \rho(r)$ for all $s\in
(r', r)$, where $r'=r-\rho(r)/2b_n$, then (2.10) follow from
(2.7).

Hence we may assume that $\rho(s)<\frac 14\rho(r)$ for some $s\in
(r', r)$. Let
$$\align
G' & =\{x\in \p B_r(0):\ \exists \ t\in (\frac sr, 1)
        \ \text{such that}\ tx\in \Ga^{\text{bdy}}_r-\Ga^{\text{bdy}}_s\},\tag 2.11\\
G''& =\{x\in \p B_r(0):\ \ \ \frac srx\in \Ga^{\text{int}}_s\}\\
  \endalign$$
be respectively the radial projection of
$\Ga^{\text{bdy}}_r-\Ga^{\text{bdy}}_s$ and $\Ga^{\text{int}}_s$ on
$\p B_r(0)$. Then $\Ga^{\text{int}}_r\subset G'\cup G''$. But since
$\rho(s)<\frac 14\rho(r)$, we have
$$\align
 |G''| & =(\frac rs)^{n-1}|\Ga^{\text{int}}_s| \\
  &\le 4^{-n+1}(\frac{r}{r'})^{n-1}|\Ga^{\text{int}}_r|\\
\endalign $$
Hence we obtain
$$|G'|\ge (1-4^{-n+1}(\frac{r}{r'})^{n-1}) |\Ga^{\text{int}}_r|.\tag 2.12$$

Regard $\Ga^{\text{bdy}}_r-\Ga^{\text{bdy}}_s$ as a (multi-valued)
radial graph over $G'$. For any point $y\in
\Ga^{\text{bdy}}_r-\Ga^{\text{bdy}}_s$, let $x$ be the projection of
$y$ on $\p B_r$. Then through the projection, the area element of
$\Ga^{\text{bdy}}_r$ at $y$ is greater than $(\frac {r'}r)^{n-1}$
times the area element of $\p B_r$ at $x$. Hence we have
$$\align
 |\Ga^{\text{bdy}}_r-\Ga^{\text{bdy}}_s| &\ge (\frac {r'}r)^{n-1}|G'|\\
           & \ge (1-4^{-n+1}(\frac{r}{r'})^{n-1})
             (\frac {r'}r)^{n-1}|\Ga^{\text{int}}_r|.\\
           \endalign $$
Note that $r'=r-\rho(r)/2b_n\to r$ as $\rho(r)\to 0$. We obtain
(2.10). $\square$

\vskip10pt

\noo{\bf Remark 2.2.} The above proof implies that if the volume
$|\Om|$ is small, then we have
$$|\Ga^{\text{bdy}}_r|\ge \frac 12 |\Ga^{\text{int}}_r|, \tag 2.13$$
where $r\in [\frac 14, 1]$. Note that in (2.13) we do not assume
that $\rho(r)$ is small. Indeed, if $\rho(r)<\min (\rho_1, \rho_2)$
is small, (2.13) is proved in Lemma 2.4. Otherwise, let
$s=r-|\Om|^{1/n}$. Define $G'$ as in (2.11) and let
$G''=\Ga^{\text{int}}_r-G'$. We have $|\Om|\ge (r-s) \big(\frac
sr\big)^{n-1}|G''| $. Hence $|G''|\le 2|\Om|^{1-1/n}$ and so
$|G'|\ge |\Ga^{\text{int}}_r|-2|\Om|^{1-1/n}$. The proof of Lemma
2.4 then implies that
$$|\Ga^{\text{bdy}}_r|\ge (\frac sr)^{n-1}|G'|
 \ge \frac 34\big(|\Ga^{\text{int}}_r|-2|\Om|^{1-1/n}\big).$$
Hence (2.13) follows if $\rho(r)\ge \min (\rho_1, \rho_2)$ and
$|\Om|<\frac 1{16} \alpha_{n-1}[\min (\rho_1, \rho_2)]^n$.

\vskip10pt

\proclaim{Lemma 2.5} Let $u$ be a smooth $H_1$-subharmonic
function in $B_1(0)$. Suppose  \newline $\sup_{B_{1/2}(0)} u>1$.
Then
$$|\{x\in B_1(0):\ \ u(x)>0\}|\ge C,\tag 2.14$$
where the constant $C>0$ depends only on $n$.
\endproclaim

\noo{\bf Proof}. We prove by contradiction, assuming that $|\Om_{1,
0}|<\delta_0^{4n}$ for some small positive constant $\delta_0$
depending on $\rho_1$ and $\rho_2$. We divide the proof into three
steps.

\noo{\bf Step 1}. For $r\in (0, 1]$ and $t\ge 0$, denote
$$\align
 & \Om_{r,t} =\{x\in B_r(0):\ \ u(x)>t\}, \tag 2.15\\
 & \Ga^{\text{bdy}}_{r,t}=\pom_{r,t}\cap B_r,\\
 & \Ga^{\text{int}}_{r,t}=\bom_{r,t}\cap \p B_r.\\
 \endalign$$
so that $\pom_{r, t}=\Ga^{\text{bdy}}_{r, t}\cup\Ga^{\text{int}}_{r,
t}$. By Sard's lemma, $\Ga^{\text{bdy}}_{r, t}$ is smooth for almost
all $t$. Note that
$$\delta_0^{4n}\ge |\Om_{1,0}|
  \ge \int_{7/8}^1|\Ga^{\text{int}}_{r, 0}| dr.\tag 2.16$$
Hence there exists $r\in [\frac 78, 1]$ such that
$|\Ga^{\text{int}}_{r, 0}|<8\delta_0^{4n}<\delta_0^{4(n-1)}$.
Without loss of generality we may also assume that
$|\Ga^{\text{int}}_{1, 0}|<\delta_0^{4(n-1)}$. Note that $\Om_{r',
t'}\subset \Om_{r, t}$ for any $r'<r, t'>t$. Hence for all $r\in (0,
1]$ and $t\ge 0$,
$$|\Om_{r, t}|<\delta_0^{4n}
 \ \ \text{and}\ \  |\Ga^{\text{int}}_{1, t}|<\delta_0^{4(n-1)}. \tag 2.17$$
To apply the previous Lemmas, we assume that $\delta_0<\frac
1{32}\min\{\rho_1, \rho_2\}$.

Consider the integration
$$0\le \int_{\Om_{r,t}} H_1[u]
 =-(\int_{\Ga^{\text{bdy}}_{r,t}}+\int_{\Ga^{\text{int}}_{r,t}})
  \frac {\gamma\cdot Du}{\sqrt{1+|Du|^2}},\tag 2.18$$
where $\gamma$ is the unit inner normal of $\Om_{r,t}$. We have
$$\align
 & |\int_{\Ga^{\text{int}}_{r,t}} \frac {\gamma\cdot Du}{\sqrt{1+|Du|^2}}|
  \le |\Ga^{\text{int}}_{r,t}|,\\
& \int_{\Ga^{\text{bdy}}_{r,t}} \frac {\gamma\cdot
Du}{\sqrt{1+|Du|^2}}
 =\int_{\Ga^{\text{bdy}}_{r,t}} \frac {|Du|}{\sqrt{1+|Du|^2}}.\\
\endalign $$
Suppose there exist $r$ and $t$ such that
$$|\Ga^{\text{bdy}}_{r,t}|\ge (1+\delta)|\Ga^{\text{int}}_{r,t}|\tag 2.19$$
for some small constant $\delta>0$ (we can fix $\delta=4^{-n}$), and
there exists a subset $\hat\Ga^{\text{bdy}}_{r,t}\subset
\Ga^{\text{bdy}}_{r,t}$ such that
$$\align
 & |Du|>2\delta^{-1/2}\ \ \ \text{on}\ \ \hat\Ga^{\text{bdy}}_{r,t},\\
 & |\hat\Ga^{\text{bdy}}_{r,t}|>(1-\frac {\delta}{4})|\Ga^{\text{bdy}}_{r,t}|.\tag 2.20\\
 \endalign $$
Then
$$\align
\int_{\Ga^{\text{bdy}}_{r,t}} \frac {|Du|}{\sqrt{1+|Du|^2}}
 & \ge \int_{\hat\Ga^{\text{bdy}}_{r,t}}
       \frac {2\delta^{-1/2}}{\sqrt{1+4\delta^{-1}}}\\
 & \ge \frac{1-\delta/4}{\sqrt{1+\delta/4}} |\Ga^{\text{bdy}}_{r,t}|\\
 & > |\Ga^{\text{int}}_{r, t}|.\\
\endalign $$
We reach a contradiction.

In the following we prove there exists $r, t$ such that (2.19) and
(2.20) hold (so we reach a contradiction and Lemma 2.5 is proved).
Accordingly we introduce the sets
$$\align
 P&=\{(r,t)\in [\frac 12, 1] \times [0, 1]:\
      |\Ga^{\text{bdy}}_{r,t}|\le (1+\delta)|\Ga^{\text{int}}_{r,t}|\},\\
 Q&=\{(r,t)\in [\frac 12, 1] \times [0, 1]:\
      |\Ga^*_{r,t}|\ge \frac {\delta}4 |\Ga^{\text{bdy}}_{r,t}|\},\\
\endalign $$
where
$$\Ga^*_{r,t}=\{x\in\Ga^{\text{bdy}}_{r,t}: \ |Du|(x)\le 2\delta^{-1/2}\}. $$
If there exists $(r, t)\in [\frac 12, 1] \times [0, 1]$ such that
$(r, t)\not\in P\cup Q$, then (2.19) and (2.20) hold and the lemma
is proved. In the following we show that both sets $P$ and $Q$
have small Lebesgue measure.

\vskip10pt

\noo{\bf Remark 2.3}. We remark that (2.19) may not hold if the
shape of $\Om_{r, t}$ is like a thumbtack, namely a flat cap with a
thin cylinder.

 \vskip10pt

\noo{\bf Step 2}. Estimate of $|Q|$. For any fixed $r\in [\frac
12, 1]$, denote
$$ Q_r=\{t\in [0, 1]:\ \ (r, t)\in Q\} $$
a slice of $Q$ at $r$, and denote $\phi(t)=|\Om_{r,t}|$. By the
co-area formula, we have, for a.e. $t$,
$$\phi'(t)
 = -\int_{\Om_{r, t}\cap\{u=t\}} \frac {1}{|Du|}
 = -\int_{\Ga^{\text{bdy}}_{r,t}} \frac {1}{|Du|},\tag 2.21$$
as $u$ is smooth and $r$ is fixed.  Hence for any $t\in Q_r$,
$$\align
\phi'(t) & \le -\int_{\Ga^*_{r,t}} \frac {1}{|Du|}\\
         & \le -\frac {1}{2\delta^{-1/2}}|\Ga^*_{r,t}|\\
         & \le -\frac {1}{8}\delta^{3/2} |\Ga^{\text{bdy}}_{r,t}|.\\
         \endalign $$
By the isoperimetric inequality,
$$|\pom_{r,t}|\ge c_n|\Om_{r,t}|^{1-\frac 1n},\tag 2.22$$
where the best constant $c_n$ is attained when the domain is a ball.
Similar to (2.16),
$$\delta_0^{4n}\ge |\Om_{1,t}|
  \ge \int_0^1|\Ga^{\text{int}}_{r, t}| dr.$$
Hence the Lebesgue measure of the set $\hat I:=\{t\in [0, 1]:\
|\Ga^{\text{int}}_{r, t}|>\delta_0^{4(n-1)}\}$ is less than
$\delta_0^{4}$. For any $t\not\in \hat I$, by Lemma 2.4,
$$|\Ga^{\text{int}}_{r,t}|\le 2|\Ga^{\text{bdy}}_{r,t}|.$$
Hence from (2.22) and noting that
$\pom_{r,t}=\Ga^{\text{bdy}}_{r,t}\cup \Ga^{\text{int}}_{r,t}$,
$$|\Ga^{\text{bdy}}_{r,t}|\ge \frac {c_n}3 |\Om_{r,t}|^{1-\frac 1n}.$$
We obtain
$$\phi'(t)\le -\frac {c_n}{24}\delta^{3/2}\phi^{1-\frac 1n}(t). \tag 2.23$$
Namely $(\phi^{\frac 1n})'(t)\le -\frac {c_n}{24n}\delta^{3/2}$ when
$t\in \hat Q_r:=Q_r-\hat I$. Note that $(\phi^{\frac 1n})'(t)\le 0$
for any $t\in [0, 1]-\hat Q_r$. Hence
$$\align
\phi^{\frac 1n}(0)-\phi^{\frac 1n}(1)
 & = -\int_{\hat Q_r}(\phi^{\frac 1n})'(t)
                -\int_{[0, 1]-\hat Q_r} (\phi^{\frac 1n})'(t)\\
 & \ge \frac {c_n}{24n}\delta^{3/2}|\hat Q_r|.\\
 \endalign  $$
We get the estimate
$$|\hat Q_r|\le \frac {24n}{c_n}\delta^{-3/2}\phi^{\frac 1n}(0).$$
By assumption,  $\phi(0)=|\Om_{1, 0}|\le \delta_0^{4n}$. Hence when
$\delta_0$ is small (recall that $\delta=4^{-n}$), we obtain $|\hat
Q_r|<\frac 12\delta_0^2$. Hence $|Q_r|\le |\hat Q_r|+|\hat I|\le
\delta_0^2$. It follows that
$$|Q|=\int_{1/2}^1|Q_r|\le \frac 12 \delta_0^2.$$
That is, $Q$ is a small set.

\vskip10pt

\noo{\bf Step 3}. Estimate of $|P|$. For any fixed $t\in [0, 1]$,
denote $P_t=\{r\in [\frac 12, 1]:\ (r, t)\in P\}$ a slice of $P$
at height $t$. We prove that $P_t$ has small Lebesgue measure, so
that $|P|=\int_0^1|P_t|$ is also small.

Denote by $\rho(r)$ the geodesic radius of $\Ga^{\text{int}}_{r,t}$,
as introduced before Lemma 2.3. Namely, we define $\rho(r)$ such
that a geodesic ball of radius $\rho(r)$ in $\p B_r$ has the volume
$|\Ga^{\text{int}}_{r,t}|$.

We first consider the case when $\rho$ is increasing in $r$. In this
case, by (2.17) we have $\rho(r)<\rho(1)<\delta_0^{4}$ for any
$r<1$. Let $\ol r_1 =\sup\, r$: $r\in [\frac 12, 1]$ and there
exists $\ul r<r$ such that
$$\align
 & \rho(\ul r)=\frac 12 \rho(r),\\
 & \frac{\rho(r)-\rho(\ul r)}{r-\ul r}\ge b_n.\tag 2.24\\
 \endalign $$
We obtain an interval $I_1=[\ul r_1, \ol r_1]$, where $\ul r_1$ is
the largest $\ul r$ satisfying (2.24). Next let $\ol r_2=\sup\, r\in
[0, \ul r_1]$ such that the above formulae hold, and we obtain an
interval $I_2=[\ul r_2, \ol r_2]$. Continue the process we obtain a
sequence of intervals $\{I_k\}$, $I_k=[\ul r_k, \ol r_k]$. By the
monotonicity of $\rho$, we have
$$\sum_k |I_k|\le \rho(1)/b_n\le \delta_0^{4 }/b_n. \tag 2.25$$

For any $r \not\in \bigcup_k I_k$ and $r\in [\frac 12, 1]$,  by
our definition of $I_k$ we have
$$\rho(s)\ge \frac 12\rho(r)\ \ \ \forall\ s\in (r', r)$$
where $r'=r-\rho(r)/2b_n$. By Lemma 2.3, $|\Ga^{\text{bdy}}_{r,
t}|\ge 2|\Ga^{\text{int}}_{r, t}|$. Hence (2.19) holds and $r\not\in
P_t$. It follows that $P_t\subset \bigcup_k I_k$. By (2.25),
$|P_t|$ is small.

\vskip10pt

Next we consider the case $\rho$ is not monotone increasing. In this
case,  we may also assume that $\sup_{r\in [\frac 14, 1]} \rho(r)$
is small. For if there exists $r_0\in [\frac 14, 1)$ such that
$\rho(r_0)\ge \delta_0^3$, we choose $r_0=\inf\{r\in [\frac 14, 1]:
\ \rho(r)\ge \delta_0^3\}$. By the argument below, the set
$P'_t=\{r\in [\frac 12, r_0]:\ (r, t)\in P\}$ is a small set. If
$r_0>1-\delta_0^2$, then $P_t\subset  P'_t\cup [r_0, 1]$ is also
small. If $r_0<1-\delta_0^2$, recall that $\delta_0^{4n}\ge
|\Om_{1,t}| \ge \int_0^1|\Ga^{\text{int}}_{r, t}| dr$. Hence the
Lebesgue measure of the set $I':=\{r\in [r_0, 1]:\
|\Ga^{\text{int}}_{r, t}|> \delta_0^{4(n-1)}\}$ is less than
$\delta_0^4$. For any $r\in [r_0, 1]-I'$, we have
$|\Ga^{\text{int}}_{r, t}|\le \delta_0^{4(n-1)}$ but by Remark 2.2,
$$|\Ga^{\text{bdy}}_{r, t}| \ge |\Ga^{\text{bdy}}_{r_0, t}|
 \ge\frac 12 |\Ga^{\text{int}}_{r_0,t}|
 \ge \frac 12\alpha_{n-1}\delta_0^{3(n-1)}$$
Hence $|\Ga^{\text{bdy}}_{r, t}| \ge 2|\Ga^{\text{int}}_{r, t}|$ and
so $r\not\in P_t$. Again $|P_t|\le |P'_t|+|I'|$ is small. In the
following we assume directly that $\sup_{r\in [1/4, 1]}
\rho(r)<\delta_0^3$.

Let
$$\hat\rho(r)=\sup\{\rho(s):\ \frac 14<s<r\},
                  \ \ \ r\in [\frac 14, 1].\tag 2.26$$
Then $\hat \rho$ is increasing in $[\frac 14, 1]$. Similarly we
define the sequence of intervals $I_k=[\ul r_k, \ol r_k]$ in terms
of $\hat\rho$. Then
$$\sum_k |I_k|\le \hat\rho(1)/b_n
   \le \delta_0^3/{b_n}.$$

We claim that $P_t\subset \bigcup_k I_k$. Indeed, for any $r
\not\in \bigcup_k I_k$ and $r\in [\frac 12, 1]$,  let
$r'=r-\rho(r)/2b_n$. If
$$\rho(s)\ge \frac 14\rho(r)\ \ \ \forall\  s\in (r', r),$$
the claim follows from Lemma 2.3. If there exists an $s\in (r',
r)$ such that $\rho(s) < \frac 14\rho(r)$, note that $r'\ge r-\hat
\rho(r)/2b_n$,  by our definition of $I_k$,
$$\hat \rho(s)\ge \frac 12 \hat\rho(r)\ge \frac 12 \rho(r)
 \ \ \ \forall\ s\in [r', r]. $$
Hence there exists $\tau \in [\frac 14, s]$  such that
$\rho(\tau)\ge \frac 12\rho(r)$. We divide $\Ga^{\text{bdy}}_{r, t}$
into three pieces, $\Ga^{\text{bdy}}_{r,
t}=\Ga_a\cup\Ga_b\cup\Ga_c$, where
$$\align
\Ga_a & =\Ga^{\text{bdy}}_{r, t}\cap\{s<|x|<r\},\\
\Ga_b & =\Ga^{\text{bdy}}_{r, t}\cap\{\tau<|x|<s\},\\
\Ga_c & =\Ga^{\text{bdy}}_{r, t}\cap\{0<|x|<\tau\}.\\
\endalign $$
By Lemma 2.4, we have
$$ \align
|\Ga_a| & \ge (1-4^{-n+5/4})|\Ga^{\text{int}}_{r, t}|,\\
|\Ga_c| & \ge (1-4^{-n+5/4})|\Ga^{\text{int}}_{\tau, t}|.\\
\endalign $$
By projecting $\Ga_b$ to $\p B_\tau$ and noticing that $\rho(s)\le
\frac 12\rho(\tau)$, we have, similarly to the proof of Lemma 2.4,
$$ |\Ga_b| \ge (1-2^{-n+5/4})|\Ga^{\text{int}}_{\tau, t}|.$$
Recall that $\sup_{r\in [\frac 14, 1]}\rho(r)\le \delta_0^3$, we
have $|\Ga^{\text{int}}_{\tau, t}|\approx
\alpha_{n+1}\rho^{n-1}(\tau)$ and $|\Ga^{\text{int}}_{r, t}|\approx
\alpha_{n+1}\rho^{n-1}(r)$. Hence by $\rho(\tau)\ge \frac 12\rho(r)$
we have $|\Ga^{\text{int}}_{\tau, t}|\ge (1-\eps)2^{-n+1}
|\Ga^{\text{int}}_{r, t}|$, where $\eps\to 0$ as $\delta_0\to 0$.
Assume $\delta_0$ small such that $\eps<4^{-n}$. Then we obtain
$$\align
|\Ga^{\text{bdy}}_{r, t}| & =|\Ga_a|+|\Ga_b|+|\Ga_c|\tag 2.27\\
 &\ge (1+4^{-n+1})|\Ga^{\text{int}}_{r, t}|.\\
 \endalign $$
The claim  $P_t\subset \bigcup_k I_k$ is proved and hence $P_t$ is a
small set. $\square$

\vskip10pt

From Lemma 2.5, we have the following weak Harnack inequality,
which is an improvement of (2.5).

\proclaim{Corollary 2.1}
 Let $u$ be an $H_1$-subharmonic function in $B_r(0)$. Then for any
constant $p>0$, there exists a constant $C$ depending on $n$ and
$p$ such that
$$\sup_{B_{r/2}} u
  \le \frac{C}{r^{n/p}} \big(\int_{B_r} (u^+)^p\big)^{1/p},\tag 2.28$$
where $u^+=\max(u, 0)$.
\endproclaim

\noo{\bf Proof}. It suffices to prove that
$$u(0)\le \frac{C}{r^{n/p}} \big(\int_{B_r} (u^+)^p\big)^{1/p}.\tag 2.29$$
We will prove it for smooth $H_1$-subharmonic functions. In the
general case it follows from the approximation in \S5.

If $\sup_{B_{r/2}} u\le r$, then (2.29) follows from (2.5). In the
following we assume that $\sup_{B_{r/2}} u\ge r$. By the
transformation $u\to u/r$ and $x\to x/r$, we may assume that $r=1$.

If $u(0)\ge 1$, applying Lemma 2.5 to the function $(u-\frac
12u(0))^+$, we see that $|\{x\in B_1(0):\ \ u(x)>\frac 12 u(0)\}|
\ge C$. Hence we obtain (2.29).

If $u(0)\le 1$, assume $\sup_{B_{1/2}} u$ is attained at $x_0$. Then
$\sup_{B_{1/2}(x_0)} u\ge 1$. Applying Lemma 2.5 to $(u-\frac
12u(x_0))^+$ in $B_{1/2}(x_0)$, we also obtain (2.29). $\square$

\vskip10pt

\vskip10pt

\noo{\bf Proof of Theorem 2.1}. Let $u$ be a nonnegative solution to
the minimal surface equation (2.3) in $B_r(0)$. It suffices to show
that $\sup_{B_{r/2}(0)} u$ is bounded from above by a constant $C$
depending only on $n, r$ and $\psi$. Once $u$ is bounded from above,
by the interior gradient estimate, equation (2.3) becomes uniformly
elliptic and the full Harnack inequality follows [GT]. Alternatively
we may also use the estimates for $\inf_{B_{1/2}(0)}u$ in [T1] or
[PS2].

By a scaling we may assume that $r=1$.  Denote $\Om_t=\{x\in
B_1(0):\ u(x)>t\}$ and $\Ga^{\text{int}}_t=\bom_t\cap\p B_1(0)$ and
$\Ga^{\text{bdy}}_t=\pom_t\cap B_1(0)$. If $\sup_{B_{1/2}(0)}u$ is
sufficiently large, by (2.14) we have $|\Om_t|\ge C$ for some $C>0$
independent of $t$. Hence by the assumption
$\lim_{t\to\infty}\psi(t)=0$, we have
$|\Ga^{\text{bdy}}_t|>2|\Ga^{\text{int}}_t|$ for all large $t$.
Namely (2.19) (with $r=1$, $\delta=1$) is satisfied for all large
$t$.

Let $\phi(t)=|\Om_t|$ and denote $Q=\{t\ge 0: \ |\Ga^*_t|\ge \frac
14 |\Ga^{\text{bdy}}_t|\}$, where $\Ga^*_t=\{x\in\Ga^{\text{bdy}}_t:
\ |Du|(x)\le 2\}$. Then from the proof of Step 2 above, $\phi$
satisfies (2.23). Hence
$$\phi^{\frac 1n}(0)-\phi^{\frac 1n}(T)
                  \ge \frac {c_n}{24n}|Q|. \tag 2.30$$
Hence (2.20) (with $r=1$, $\delta=1$) is satisfies for most large
$t$. Choosing a $t\not\in Q$, we reach a contradiction as in Step
1 of the proof of Lemma 2.5. $\square$

\vskip10pt

\noo{\bf Remark 2.4}. From the proof of Lemma 2.5 (see (2.18)), one
sees that if for any $\omega\subset \Om$,
$$\int_{\omega} H_1[u] \ge -\nu(\omega)\tag 2.31$$
for some nonnegative measure $\nu$ satisfying $\frac
{\nu(\omega)}{|\p\omega|}\to 0$ as $|\omega|\to 0$, then estimate
(2.14) holds, with the constant $C$ depending also on $\nu$. This
estimate, combined with Theorem 3.1 in [T1],  implies a Harnack
inequality for solutions $u\in W^{2, n}(\Om)$ to the non-homogeneous
mean curvature equation.

 \vskip20pt

\centerline{\bf 3. Gradient and uniform estimates}

 \vskip10pt

First we establish an integral gradient estimate.

\proclaim{Theorem 3.1} Let $u\in C^2(\Om)$ be a non-positive
$H_1$-subharmonic function. Then for any open set $\omega\Subset
\Om$,
$$\int_{\omega}|D u_t|\leq C,\tag 3.1$$
where $u_t=\max(u, -t)$, $t$ is a constant, and $C>0$ depends on
$\omega, t$, but is independent of $u$.
\endproclaim

\noo{\bf Proof}. Let $\phi(x)\in C_0^{\infty}(\Om)$ be a smooth
function with support in $\Om$ such that $0\leq\phi(x)\leq 1$ and
$\phi(x)\equiv 1$ on $\omega$. We may assume that $|\pom|$, the
area of $\pom$, is bounded, otherwise we may restrict to a
subdomain of $\Om$ which contains $\omega$. Then
$$ \align
 \int_{\Om}\phi (-u_t)  H_1[u]
 &= \int_{\Om}\frac{\phi|D u_t |^2}{\sqrt{1+|D u_t |^2}}+
    \int_{\Om}\frac{ u_t D u_t \cdot D\phi}{\sqrt{1+|D u_t |^2}}\\
 &\geq\int_{\omega}\frac{|D u_t |^2}{\sqrt{1+|D u_t |^2}}+
  \int_{\Om}\frac{ u_t D u_t \cdot D\phi}{\sqrt{1+|D u_t |^2}}\\
 &\geq\int_{\omega}|D u_t |-|\omega|
  +\int_{\Om}\frac{ u_t D u_t  \cdot D\phi} {\sqrt{1+|D u_t |^2}}.
\endalign$$
Note that
$$\int_{\Om}\phi (-u_t) H_1[u]\leq t\int_\Om H_1[u]\le t|\pom|$$
and
$$ \int_{\Om}\frac{ u_t D u_t \cdot
 D\phi}{\sqrt{1+|D u_t |^2}} \leq Ct|\Om|.$$
We obtain
$$\int_{\omega}|D u_t |\le C(1+t)(|\Om|+|\pom|). $$
Hence (3.1) is proved. $\square$

In the next section we will prove that every $H_1$-subharmonic
function can be approximated by smooth ones. Note that if $u\in
SH_1(\Om)$, then $u_t\in SH_1(\Om)$.  Hence by Theorem 3.1 we have

\proclaim{Corollary 3.1} For any $u\in SH_1(\Om)$ bounded from
above and any $\Om'\Subset\Om$, $u_t\in BV(\Om')$. In particular
if $u$ is bounded from below, then $u\in BV(\Om')$.
\endproclaim

By the example in \S 8,  $u\not\in W^{1,1}(\Om')$ in general.

Next we consider the $L^\infty$ estimate for $H_1$-subharmonic
functions. We say a set $A$ is Caccioppoli if it is a Borel set with
characteristic function $\phi_A$ whose distributional derivatives
$D\phi_A$ are Radon measures [G3]. If $A$ is Caccioppoli, we have
$$|\p A|=\int_{R^n} |D\phi_A|. \tag 3.2$$

\proclaim{Theorem 3.2} Assume that $u\in SH_1(\Om)\cap C^2(\Om)$
is bounded from below on $\pom$. Assume that there is a positive
constant $\eta$ such that for any Caccioppoli set $A\subset\Om$,
$$\int_A H_1[u]\leq(1-\eta)|\p  A|.\tag 3.3$$
Then there is a constant $C>0$ such that
$$\inf\limits_{x\in\Om}u\geq -C.\tag 3.4$$
\endproclaim

\noo {\bf Proof}. For any $t>0$, denote $\Om_t=\{x\in\Om: \
u(x)\leq -t\}$ and $\p_1\Om_t=\{x\in\p \Om_t: \ \ |Du|\leq
t^{2/3}\}$. Since $u$ is bounded from below on $\p \Om$, we may
choose a large $T$ such that $\Om_T\Subset\Om$ and
$$\frac{T^{2/3}}{\sqrt{1+T^{4/3}}}\geq 1-\eta/2.\tag 3.5$$
We claim that for any $t>T$,
$$|\p_1\Om_t|\geq\frac \eta 2 |\p \Om_t|.\tag 3.6$$
Indeed, if there exists a $t\ge T$ such that
$|\p_1\Om_{t}|<\frac{\eta}{2}|\p \Om_{t}|$, we have
$$\align
\int_{\Om_{t}}H_1[u]
 &=\int_{\p \Om_{t}} \frac{|Du|}{\sqrt{1+|Du|^2}}\\
 &\geq\int_{\p \Om_{t}-\p_1\Om_{t}}
                             \frac{|Du|}{\sqrt{1+|Du|^2}}\\
 &\geq(1-\eta/2)(1-\eta/2)|\p \Om_{t}|\\
 &>(1-\eta)|\p \Om_{t}|,
\endalign $$
which is in contradiction with the assumption (3.3).

Let $\varphi(t)=|\Om_t|$. If $t>-\inf_{\pom} u$, then
$\Om_t\subset\subset \Om$. Hence by the co-area formula,
$$\varphi'(t)=-\int_{\p \Om_t}\frac{1}{|Du|}
 \le -\int_{\p_1\Om_t}\frac{1}{|Du|}.$$
When $t>T$,
$$\varphi'(t)\leq -\frac{\eta}{2t^{2/3}}|\p \Om_t|.$$
By the isoperimetric inequality,
$$\varphi^{1-1/n}(t)\leq C|\p \Om_t|,$$
we obtain
$$\varphi'(t)\leq -\frac{C\eta}{t^{2/3}}\varphi^{1-1/n}(t).\tag 3.7$$
Namely $[\varphi^{1/n}(t)]'\leq -C\eta t^{-2/3}$. Taking
integration from $T$ to $t$, we obtain
$$\varphi^{1/n}(t)\leq\varphi^{1/n}(T)
       +C\eta (T^{\frac{1}{3}}-t^{\frac{1}{3}})\tag 3.8$$
for a different $C$.  Hence $\varphi$ vanishes when
$t>C\big[T+\big(\frac{|\Om|^{1/n}}{\eta}\big)^3\big]$. This
completes the proof. $\square$

\vskip10pt

\noo{\bf Remark 3.1}. Condition (3.3) was introduced in [Gia], in
which it is proved that (1.3) is necessary and sufficient for the
existence of a minimizer of an associated functional.  From the
proof of Theorem 3.2 one sees that the condition (3.3) can be
weakened to
$$|G(\ol t)|\to \infty \ \ \ \text{as}\ \ \ol t\to \infty,\tag 3.9$$
where $G(\ol t)$ is the set of $t\in (0, \ol t)$ such that
$$\int_{\Om_t} H_1[u]\leq (1-\eta)|\pom_t|.\tag 3.10$$
This is because (3.6) and (3.7) hold for any $t\in G(\ol t)$.
Furthermore, as the co-area formula holds for BV functions [G3],
the above argument applies to BV functions.

\noo{\bf Remark 3.2}. From the proof, the constant $C$ in Theorem
3.2 depends only on $n, \Om$, $\eta$, and $\inf_{\pom} u$. Hence
Theorem 3.2 also holds for non-smooth $H_1$-subharmonic functions,
by the approximation in Section 5.

\noo{\bf Remark 3.3}. A similar estimate for the prescribing
$k$-curvature equation was established in [T2]. We include a direct
proof for the mean curvature case (namely the case $k=1$) here for
completeness.

\vskip20pt

\centerline{\bf 4. Perron lifting }

\vskip10pt

Let $u$ be an $H_1$-subharmonic function in $\Om$ and let
$\omega\Subset\Om$ be an open, precompact set in $\Om$.  The
Perron lifting of $u$ in $\omega$, $u^{\omega}$, is defined as the
upper semicontinuous regularization of
$$ u^*=\sup\{v\ |\ v \ \text{is $H_1$-subharmonic in $\Om$ and}
 \ v\le u\ \text{in}\ \ \Om-\omega\}, \tag 4.1$$
namely
$$u^\omega(x)=\lim\limits_{r\rightarrow 0}\sup\limits_{B_r(x)}u^*. \tag 4.2$$

\noo{\bf Remark 4.1}.  Obviously we have $u^\omega\ge u$ on
$\p\omega$. However for general open set $\omega$, it may occur
that $u^\omega>u$ on part of the boundary $\p\omega$, even if $u$
is a smooth function. This is easily seen by considering the
Perron lifting in $\omega=B_{R}-\ol B_r$ of a radial function $u$,
where $R>r$. Then in general one has $u^\omega>, \ne u$ on the
inner boundary $\p B_r$. But if $u$ is continuous, by Lemma 2.2
one has $u^\omega=u$ on the outer boundary $\p B_{R}$.

First we prove the following basic result for $H_1$-harmonic
functions. Note that our definition of $H_1$-harmonic functions
does not imply they are bounded from below.

\proclaim{Lemma 4.1}
 Let $u$ be an $H_1$-harmonic function in $\Om$.
Then $u$ is locally bounded and smooth in $\Om$, and satisfies the
equation $H_1[u]=0$ in $\Om$.
\endproclaim

\noo{\bf Proof}. Assume that $B_1(0)\Subset \Om$. By definition,
an $H_1$-harmonic function is $H_1$-subharmonic.  The
$n$-dimensional Hausdorff measure $|\{x\in \Om:\ u<-t\}|\to 0$ as
$t\to\infty$. Hence we may assume that the $(n-1)$-dimensional
Hausdorff measure $|\{x\in \p B_1 :\ u<-t\}|\to 0$ as
$t\to\infty$.

Since $u$ is upper semicontinuous, there exists a sequence of
smooth functions $\{v_j\}$ in $\Om$ such that $v_j\searrow u$,
namely $v_j$ converges to $u$ monotone decreasingly. By Lemma 2.2,
there is a solution $\hat v_j\in C^2(B_1)\cap C^0(\ol B_1)$ to
$$\cases
  H_1[v]=0\ \ \ &\text{in}\  B_1(0),\\
 v=v_j &\text{on}\ \ \p B_1.
 \endcases \tag 4.3$$
Since $\hat v_j$ is monotone decreasing and $\hat v_j>u$, it is
convergent. We may assume that $\hat v_j\searrow \hat v$.
Obviously $\hat v\ge u$ in $B_1$.

Next we show that $\hat v\le u$ on $\p B_1$, namely for any given
$x_0\in\p B_1$,
$$\lim_{x\to x_0} \hat v(x)\le u(x_0), \tag 4.4$$
so that
$$\hat v\equiv u\ \ \ \text{in}\ \ B_1. $$
 Indeed, since $u$ is upper semicontinuous on $\p B_1$,
there is a continuous function $w$ on $\p B_1$ such that
$w(x_0)=u(x_0)$ and $w\ge u$ on $\p B_1$. By the monotonicity of
$v_j$ on $\p B_1$, it is easy to show that for any $\eps>0$, there
is a $\delta>0$ such that for sufficiently large $j$,
$v_j(x)<u(x)+\eps$ in $\{x\in\p B_1:\ |x-x_0|\le \delta\}$. Hence
by adding $C|x-x_0|^2$ to $w$ for some large $C$,  we may assume
that $w>v_j-\eps$ on $\p B_1$ when $j$ is sufficiently large. Let
$\hat w\in C^2(B_1)\cap C^0(\ol B_1)$ be the solution of $H_1[\hat
w]=0$ in $B_1(0)$, satisfying the boundary condition $\hat w=w$ on
$\p B_1$. Then $\hat w\ge \hat v_j-\eps\ge \hat v-\eps$. Since
$\eps>0$ is arbitrary, we obtain $u(x_0)=\hat w(x_0)\ge \hat
v(x_0)$, namely (4.4) holds.

If $\inf_{B_{1/2}} \hat v_j\to-\infty$ as $j\to\infty$, by the
Harnack inequality (Theorem 2.1), we see that $\hat v_j\to-\infty$
uniformly in $B_{1/2}$. Recall that $\hat v_j\ge u$. We obtain
$u=-\infty$ in $B_{1/2}$. But by the definition of subharmonic
functions, the set $\{u=-\infty\}$ has measure zero. We reach a
contradiction. Hence $\hat v_j$ is locally uniformly bounded, and so
$u$ is locally uniformly bounded and smooth. Note that to apply
Theorem 2.1 we need the condition $|\{u(x)<-t:\  x\in \p
B_1\}|_{n-1}\to 0$ as $t\to\infty$, which is satisfied as noted at
the beginning of the proof. $\square$

\vskip10pt

\noo{\bf Remark 4.2}. The function $\hat v$ is independent of the
sequence $v_j$. Indeed, let $w_j$ be another sequence of smooth
functions on $\p B_1$ such that $w_j\searrow u$. Let $\hat w_j$ be
the solution of (4.3) with boundary condition $\hat w_j=w_j$ on
$\p B_1$ and let $\hat w=\lim \hat w_j$. Then by (4.4), we have
$\hat w_j\ge \hat v$. Hence $\hat w\ge \hat v$. Similarly we have
$\hat v\ge \hat w$. Therefore we may regard $\hat v$ as the
solution of the Dirichlet problem $H_1[v]=0$ in $B_1$ with $v=u$
on $\p B_1$.

\proclaim{Lemma 4.2}
 Let $u\in SH_1(\Om)$. Then for any open set $\omega\Subset\Om$,
 the Perron lifting $u^{\omega}$ is $H_1$-harmonic in $\omega$
 and $H_1$-subharmonic in $\Om$.
\endproclaim

\noo {\bf Proof}. The property that $u^\omega$ is
$H_1$-subharmonic in $\Om$ follows by definition. Indeed, let
$E\subset \Om$ be an open set and $h\in C^2(\ol E)$ be an
$H_1$-harmonic function satisfying $h\ge u^\omega$ on $\p E$. Then
for any $H_1$-subharmonic function $v$ in (4.1),  $h\ge v$ on $\p
E$. Hence $h\ge v$ in $E$. By the definition of $u^\omega$ in
(4.1) and (4.2) and note that $h\in C^2(E)$, it follows that $h\ge
u^\omega$ in $E$. That is, $u^\omega$ is $H_1$-subharmonic.

To show that $u^\omega$ is $H_1$-harmonic in $\omega$, let
$B_r\Subset \omega$ and let $v$ be the solution of the Dirichlet
problem $H_1[v]=0$ in $B_r$ with $v=u^\omega$ on $\p B_r$ (see
Remark 4.2). Then $v\ge u^\omega$ in $B_r$. Let $\hat u=v$ in
$B_r$ and $\hat u=u^\omega$ in $\Om-B_r$. Then $\hat u$ is upper
semicontinuous and $H_1$-subharmonic. It follows by (4.1) that
$\hat u\le u^\omega$. Hence $u^\omega=v$ in $B_r$. Namely
$u^\omega$ is $H_1$-harmonic in $B_r$. $\square$

\proclaim{Lemma 4.3} Suppose $\{u_j\}\subset SH_1(\Om)$ such that
$u_j$ converges to a measurable function $u$ a.e. with
$|\{u=-\infty\}|=0$. Let $\wtt u$ be the upper semicontinuous
regularization of $u$. Then $\wtt u=u$ a.e. and $\wtt u$ is
$H_1$-subharmonic.
\endproclaim

\noo {\bf Proof}. Let $x_0$ be a Lebegue point of $u$. By adding a
constant we assume that $u(x_0)=0$. Then Lemma 2.5 implies that
$\sup_{B_r(x_0)} u\to 0$ as $r\to 0$. Hence $u=\wtt u$ at all
Lebegue points, namely $u=\wtt u$ a.e..

To prove that $\wtt u$ is $H_1$-subharmonic, let $\omega\Subset
\Om$ be an open set and $h\in C^2(\ol\omega)$ be an $H_1$-harmonic
function with $h\ge \wtt u$ on $\p\omega$. If $u_j$ is monotone
decreasing, then for any $\eps>0$, by the monotonicity and the
upper semicontinuity of $u_j$, $h\ge u_j-\eps$ on $\p\omega$
provided $j$ is sufficiently large. It follows that $h\ge
u_j-\eps$ in $\omega$ for all large $j$. Hence $h\ge \wtt u$ in
$\omega$ and so $\wtt u$ is $H_1$-subharmonic. If $u_j$ is
monotone increasing, obviously $h\ge u_j$ on $\p\omega$ for all
$j$. Hence $h\ge \wtt u$ in $\omega$ and so $\wtt u$ is
$H_1$-subharmonic.

For general $\{u_j\}$, let $w_{k, j}=\max\{u_k, \cdots, u_j\}$.
Then for fixed $k$, $w_{k, j}\nearrow w_k$ a.e., as $j\to\infty$,
for some $w_k\in SH_1(\Om)$, and $w_k\searrow u$ a.e. as
$k\to\infty$. Hence $u$ is $H_1$-subharmonic. $\square$

For $u\in SH_1(\Om)$, the Perron lifting $u^{B_t}$ is monotone
increasing in $t$,
$$ \lim\limits_{t\rightarrow\delta^-}u^{B_t}\leq
  u^{B_\delta}(x)\leq\lim\limits_{t\rightarrow\delta^+}u^{B_t}\
  \ \forall \ x\in\Om. \tag 4.5$$
This implies that $\|u^{B_t}\|_{L^1(\Om)}$, as a function of $t$,
is monotone and bounded. Hence, $\|u^{B_t}\|_{L^1(\Om)}$ is
continuous for almost all $t$. Since $u^{B_t}$ is continuous in
$B_t$, it follows that
$$ \lim\limits_{t\rightarrow r}u^{B_t}(x)=u^{B_r}(x)\
  \ \text{for a.e.}\ \ r>0. \tag 4.6 $$
Similar to Lemma 3.6 in [TW4], we have the following

\proclaim{Lemma 4.4} Suppose $u_j,\ u\in SH_1(\Om)$ and
$u_j\rightarrow u$ a.e. in $\Om$. Then for any $B_r\Subset\Om$
such that (4.6) holds, we have $u^{B_r}_j\rightarrow u^{B_r}$ a.e.
in $\Om$ as $j\rightarrow\infty$.
\endproclaim

\noo {\bf Proof}.\ Since $u^{B_r}_j$ and $u^{B_r}$ are locally
uniformly bounded in $C^2_{loc}(B_r)$, by passing to a
subsequence, we may assume that $u^{B_r}_j$ is convergent. Let
$w'=\lim u^{B_r}_j$ and $w$ be the upper semicontinuous
regularization of $w'$ (note that $w$ and $w'$ can differ only on
$\p B_r$). Then $w\in SH_1(\Om)$ and $w=u$ in $\Om-\ol B_r$. Hence
by the definition of the Perron lifting, we have $u^{B_r}\geq w$.

Next we prove that for any $\delta>0$,  $w\ge u^{B_{r-\delta}}$.
Once this is proved, we have $u^{B_r}\ge w\ge  u^{B_{r-\delta}}$.
Sending $\delta\to 0$, we obtain $u^{B_r}= w$ by (4.6).

To prove $w\ge u^{B_{r-\delta}}$, it suffices to prove that for
any $\eps>0$, $u^{B_r}_j\ge u-\eps$ on $\p B_{r-\delta}$ for
sufficiently large $j$.  By the interior gradient estimate,
$u^{B_r}_j$ is uniformly bounded in $C^2( B_{r-\delta/4})$. If
there exists a point $x_0\in \p B_{r-\delta}$ such that
$u(x_0)>u^{B_r}_j(x_0)+\eps$ for all large $j$, by Lemma 2.5,
there is a Lebesgue point $x_1\in B_{\delta/4}(0)$ of $u$ such
that $u(x_1)>u^{B_r}_j(x_1)+\frac 12\eps$ for all large $j$. It
follows that the limit function $w=\lim_{j\to\infty} u^{B_r}_j$ is
strictly less than $u$ a.e. near $x_1$. We reach a contradiction
as $w=\lim_{j\to \infty}u^{B_r}_j\ge \lim_{j\to\infty} u_j=u$.
$\square$

\vskip20pt

\centerline{\bf 5. Approximation by smooth functions}

\vskip10pt

We prove that every $H_1$-subharmonic function can be approximated
by a sequence of smooth, $H_1$-subharmonic functions.

\proclaim{Theorem 5.1} For any $u\in SH_1(\Om)$, there is a
sequence of smooth functions $\{u_j\}\subset SH_1(\Om)$ such that
$u_j\rightarrow u$ a.e. on $\Om$.
\endproclaim

\noo{\bf Proof}. For each $j=1, 2, \cdots$, let $\{B_{j, k}, k=1,
2, \cdots, k_j\}$ be a family of finitely many balls of radius
$2^{-j}$, contained in $\bom$, such that $\Om_{2^{-j-1}}\subset
\cup_{k=1}^{k_j} B_{j, k}$, where $\Om_\delta =\{x\in\Om:\ \
\text{dist}(x, \pom)>\delta\}$.

Let $u_{j, 0}=u$. For $m=1, \cdots, k_j$, define $u_{j, m}$ such
that $u_{j, m}=u_{j, m-1}$ in $\Om-B_{j, m}$ and $u_{j, m}$ is the
solution of
$$\cases
 H_1[v]=0\ \ \ & \text{in}\ \ B_{j, m},\\
 v=u_{j, m-1}\ \ &\text{on}\ \p B_{j, m}.\\
 \endcases \tag 5.1$$
and denote $u_j=u_{j, k_j}$. Then $u_j$ is a sequence of piecewise
smooth $H_1$-subharmonic functions and
$$u_j\ge u. \tag 5.2$$

To show that $u_j\to u$ a.e., recall that every upper
semi-continuous function $u$ can be approximated by a sequence of
smooth, monotone decreasing functions $\{v^m\}$, namely
$v^m\searrow u$. For each $m$, define $v^m_j$ as above.  Then we
have $v^m_j\to v^m$ as $j\to\infty$. Hence we may choose $j=j_m$
large such that $v^m_{j_m}\to u$ a.e.. Note that $v^m_{j_m}\ge
u_{j_m}$. Hence $u_j\to u$ a.e..

In the above proof we obtain a sequence of piece-wise smooth
functions $\{u_j\}\subset SH_1(\Om)$ which converges to $u$. To
prove the theorem we make certain mollification of $u_{j, k}$. A
simple way is to replace $u_{j, k}$ by the convolution $u_{j,
k}*\rho_{\eps}$ ($\eps$ depends on $j, k$, and $\eps\to 0$
sufficiently fast as $j\to\infty$), where
$\rho_\eps=\eps^{-n}\rho(\frac x\eps)$ and $\rho$ is a mollifier.
Namely $\rho$ is a nonnegative function satisfying $\rho \in
C^\infty_0(B_1(0))$ and $\int_{B_1}\rho=1$. Specifically we may
choose
$$\rho(x)=\cases
C exp(\frac{1}{|x|^2-1}) &\text{for}\ \ |x|\leq 1 ,\\
0 & \text{for}\ \ |x|\geq 1,
\endcases \tag 5.3$$
where $C$ is chosen such that $\int_{R^n}\rho(x)dx=1$.

The function $u_{j, k}*\rho_{\eps}$ may not be $H_1$-subharmonic.
But we have
$$H_1[u_{j, k}*\rho_{\eps}]\ge -\delta\tag 5.4$$
with $\delta\to 0$ as $\eps\to 0$. This is fine for our treatment,
as the mean curvature operator is elliptic for any smooth
functions.

We can also mollify $u_{j, k}$ in the following way to get a
sequence of $C^{1,1}$ smooth, $H_1$-subharmonic functions which
converges to $u$. For a fixed $j$, recall that we first get the
function $u_{j, 1}$, which is smooth in $B_{j, 1}$. We then get
$u_{j, 2}$, which is the Perron lifting of $u_{j, 1}$ in $B_{j,
2}$. The function $u_{j, 2}$ is piece-wise smooth in $B_{j,1}\cup
B_{j,2}$, its gradient may have a jump across the boundary
$\Ga=:B_{j, 1}\cap\p B_{j, 2}$.  If $Du_{j, 2}$ has a jump at some
point on $\Ga$, then by the maximum principle, we have $u_{j,
2}>u_{j, 1}$ in $B_{j, 2}-B_{j, 1}$.  By the Hopf lemma, $Du_{j,
2}$ has a jump at every point on $\Ga$.

Let us indicate the mollification of $u_{j, 2}$ near $\Ga$.  By a
proper choice of the axes, we assume that $B_{j, 2}$ is centered
at $(0, 2^{-j})$ and $B_{j, 1}$ is centered at $(0, c)$ for some
$c<2^{-j}$. Then $\Ga$ is given by
$$x_n= g(x')=2^{-j}-\sqrt{2^{-2j}-|x'|^2},\tag 5.5 $$
where $x'=(x_1, \cdots, x_{n-1})$. Let
$$ \align
 a(x')
  &=\lim_{t\to +0}
  \frac 1t\big[u_{j, 2}(x', g(x')+t)-u_{j, 2}(x', g(x'))\big] \tag 5.6\\
  & \ \ \ -\lim_{t\to +0}
  \frac 1t\big[u_{j, 2}(x', g(x'))-u_{j, 2}(x', g(x')-t)\big]\\
  &=\p_{x_n} u_{j, 2}(x', g(x'))-\p_{x_n} u_{j, 1}(x', g(x')) .\\
\endalign $$
By the Hopf lemma, $a(x')>0$ for all $x'$ near $0$. Let
$$\phi(x)=\frac {a(x')}{4\eps}(x_n-g(x')+\eps)^2,\tag 5.7$$
where $\eps<<2^{-j}$ is a small constant. Now let
$$\wtt u_{j, 2}(x)=
\cases
 u_{j, 2}(x) &\text{if}\ \   |x_n-g(x')|\ge \eps,\\
 u_{j, 2}(x)+\phi(x) & \text{if}\ \
                      g(x')-\eps\le x_n\le g(x'),\\
 u_{j, 2}(x)+\phi(x)-a(x')(x_n-g(x')) & \text{if}\ \
                      g(x')\le x_n\le g(x')+\eps,\\
  \endcases\tag 5.8
$$
It is obvious that $\wtt u_{j, 2}\in C^{1,1}$.  When
$g(x')-\eps\le x_n\le g(x')+\eps$,
$$D^2 \phi=\frac {a}{2\eps}(Dg, -1)\otimes(Dg,-1)+O(1).\tag 5.9$$
Note that
$$H_1[u]=\text{trace of }
  \big(1-\frac{u_iu_j}{1+|Du|^2}\big)\big( D^2 u\big)\tag 5.10$$
and the matrix $\big(1-\frac{u_iu_j}{1+|Du|^2}\big)$ is positive
definite (since $|Du|\le C$).  Hence $\wtt u_{j, 2}$ is
$SH_1$-subharmonic when $\eps$ is sufficiently small.

After the modification, $u_{j, 2}$ is smooth in $B_{j,1}\cup B_{j,
2}$. Next we can modify $u_{j, k}$, for $k=3, 4, \cdots$, in the
same way, but the constant $\eps$ will be chosen smaller and
smaller. $\square$

We note that by choosing the function $\phi$ in (5.7) more
carefully, one can make the function $\wtt u_{j, 2}$ in (5.8)
$C^{2, 1}$-smooth.

\vskip20pt

\centerline{\bf 6.  Weak convergence }

\vskip5pt

For $u\in SH_1\cap C^2$, denote $\mu_1[u]=H_1[u]dx$ the associated
measure. In this section, we prove the following weak convergence
result for $H_1[u]$.

\proclaim{Lemma 6.1} Let $u_j\in C^2(\Om)$ be a sequence of
$H_1$-subharmonic functions which converges to $u\in SH_1(\Om)$
a.e. in $\Om$. Then $\{\mu_1[u_j]\}$ converges to a measure $\mu$
weakly.
\endproclaim

\noo{\bf Proof}. For any open set $\omega\subset\Om$,
$$\mu_1[u_j](\omega)\le \mu_1[u_j](\Om)\le |\pom|\tag 6.1$$
is uniformly bounded. Hence there is a subsequence of $\mu_1[u_j]$
which converges weakly to a measure $\mu$. We need to prove that
$\mu$ is independent of the choice of subsequences of $\{u_j\}$.

Let $\{u_j\}, \{v_j\}\subset SH_1(\Om) \bigcap C^2(\Om)$. Suppose
both sequences converge to $u$ a.e. in $\Om$ and
$$\mu_1[u_j]\rightarrow\mu,\ \ \ \ \ \
  \mu_1[v_j]\rightarrow\nu \tag 6.2$$
weakly as measures. We claim that for any ball $B_r(x_0)$ such
that $B_{2r}(x_0)\subset\Om$,
$$\mu(B_r)=\nu(B_r),\tag 6.3$$
or equivalently, for any $t>0$,
$$
\align
& \mu(B_r) \leq\nu(B_{r+t}),\tag 6.4a\\
& \nu(B_r)\leq\mu(B_{r+t}).\tag 6.4b\\
\endalign
$$

We choose finitely many small balls $\{B_l\}_{l=1}^k$ contained in
$B_{r+4t/5}-B_{r+t/5}$ such that the center of each ball is on $\p
B_{r+t/2}$ and $\ol {B}_{r+3t/4}-B_{r+t/4} \subset
\bigcup\limits_{l=1}^kB_l$. Now let $u_{j,1}$ be the Perron
lifting of $u_j$ on $B_1$, and let $u_{j,2}$ be the Perron lifting
of $u_{j,1}$ on $B_2$, $\cdots$, and let $u_{j,k}$ be the Perron
lifting of $u_{j,k-1}$ on $B_k$. Denote $u_j^t=u_{j, k}$.
Similarly we obtain $v^t_j$ and $u^t$. Then $u^t_j, v^t_j$ and
$u^t$ are piece-wise smooth in $B_{r+3t/4}-B_{r+t/4}$, and
$u^t_j=u_j$, $v^t_j=v_j$ in $B_r$, and so are smooth in $B_r$. By
Lemma 4.4, we have
$$u^t_j, v^t_j\rightarrow u^t\ \ \ \ \text{in}\ \ \ \Om\ a.e.\tag 6.5$$
and
$$Du^t_j, Dv^t_j\to Du^t\ \ \ \text{on} \ \ \ \p B_{r+t/2}\ a.e.$$

Let $u^t_{j, \eps}=u^t_j*\rho_\eps$ and $v^t_{j,
\eps}=v^t_j*\rho_\eps$ be the mollifications of $u^t_j$ and
$v^t_j$, where $\rho_\eps=\eps^{-n}\rho(\frac x\eps)$ and $\rho$
is a mollifier, as was given in (5.3). Then $H_1[u^t_{j, \eps}]\ge
-\delta_\eps$ with $\delta_\eps\to 0$ as $\eps\to 0$. Noting that
$u^t_j$ is independent of $t$ in $B_r$, we have
$$
\align
  \int_{B_r} H_1[u_j]=\int_{B_r} H_1[u^t_j]
 &=\lim_{\eps\to 0} \int_{B_r} H_1[u^t_{j, \eps}] \tag 6.6\\
 &\le \lim_{\eps\to 0} \int_{B_{r+t/2}} H_1[u^t_{j, \eps}]\\
 &=\lim_{\eps\to 0} \int_{\p B_{r+t/2}} \frac
  {\gamma\cdot Du^t_{j, \eps}}{\sqrt {1+|Du^t_{j, \eps}|^2}}\\
 &= \int_{\p B_{r+t/2}}
  \frac{\gamma\cdot Du^t_{j}}{\sqrt {1+|Du^t_{j}|^2}},\\
  \endalign
$$
where $\gamma$ denotes the unit outer normal. Recall that $u^t_j,
v^t_j$ and $u^t$ are piece-wise smooth in $B_{r+3t/4}-B_{r+t/4}$,
we have
$$\lim_{j\to\infty} \int_{\p B_{r+t/2}}
  \frac{\gamma\cdot Du^t_{j}}{\sqrt {1+|Du^t_{j}|^2}}
  = \int_{\p B_{r+t/2}}
  \frac{\gamma\cdot Du^t}{\sqrt {1+|Du^t|^2}}. \tag 6.7$$
Similarly we have
$$\align
 \int_{\p B_{r+t/2}}
  \frac{\gamma\cdot Du^t}{\sqrt {1+|Du^t|^2}}
  &=\lim_{j\to\infty} \int_{\p B_{r+t/2}}
  \frac{\gamma\cdot Dv^t_{j}}{\sqrt {1+|Dv^t_{j}|^2}} \tag 6.8\\
  &=\lim_{j\to\infty} \lim_{\eps\to 0} \int_{\p B_{r+t/2}}
   \frac{\gamma\cdot Dv^t_{j, \eps}}{\sqrt {1+|Dv^t_{j, \eps}|^2}}.\\
\endalign $$
Note that
$$ \int_{\p B_{r+t}} \frac
  {\gamma\cdot Dv^t_{j, \eps}}{\sqrt {1+|Dv^t_{j, \eps}|^2}}
 - \int_{\p B_{r+t/2}} \frac
  {\gamma\cdot Dv^t_{j, \eps}}{\sqrt {1+|Dv^t_{j, \eps}|^2}}
 =\int_{B_{r+t}-B_{r+t/2}} H_1[v^t_{j, \eps}]$$
and $H_1[v^t_{j, \eps}]\ge -\delta_\eps$ with $\delta_\eps\to 0$
as $\eps\to 0$. Hence the right hand side of (6.8)
$$
  \le \lim_{j\to\infty} \lim_{\eps\to 0} \int_{\p B_{r+t}}
  \frac{\gamma\cdot Dv^t_{j, \eps}}{\sqrt {1+|Dv^t_{j,\eps}|^2}}.$$
Note that $v^t_{j,\eps}$ is independent of $t$ on $\p B_{r+t}$.
The above formula
$$ =\lim_{j\to\infty} \lim_{\eps\to 0}
              \int_{ B_{r+t}} H_1[v_{j,\eps}]
  =\lim_{j\to\infty} \int_{ B_{r+t}} H_1[v_{j}].$$
Hence we obtain $ \mu(B_r)\leq\nu(B_{r+t})$. Similarly, we can
prove $\nu(B_r)\leq\mu(B_{r+t})$. This completes the proof.
$\square$

From the above lemma, we can assign a measure $\mu$ to $u$ for any
$u\in SH_1(\Om)$, and obtain the following weak convergence
theorem.

\proclaim{Theorem 6.1}
 For any $u\in SH_1(\Om)$, there exists a Radon measure
$\mu_1[u]$ such that
\newline
(i) $\mu_1[u]=H_1[u]dx$ if $u\in C^2(\Om)$,
\newline
(ii) if $\{u_j\}\subset SH_1(\Om)$ is a sequence which converges
to $u$ a.e., then $\mu_1[u_j]\rightarrow\mu_1[u]$ weakly as
measure.
\endproclaim

Note that in (ii) above, we need to use the approximation in
Section 5.

\noo{\bf Remark 6.1}. If $\{u_j\}$ is a sequence of semi-convex
functions converging to $u$, then the weak convergence $\mu_1[u_j]
\rightharpoonup \mu_1[u]$ is a special case of the weak continuity
of Federer [F1].

\vskip20pt

\centerline{\bf 7. Existence of weak solution}

\vskip10pt

In this section we consider the Dirichlet problem
$$\align
H_1[u]& =\nu\ \ \ \text{in}\ \ \Om,\tag 7.1\\
u & =\varphi\ \ \text{on}\ \ \p \Om,\\
\endalign$$
where $\Om$ is a bounded smooth domain in $\R^n$, $\phi$ is a
continuous function on $\pom$, and $\nu$ is a nonnegative measure.
Here we also use $\nu$ to denote its density with respect to the
Lebesgue measure.

For the Dirichlet problem of the mean curvature equation, usually
one assumes that the right hand side $\nu$ is Lipschitz continuous
so that the solution is smooth [GT, G1]. When $\nu\in L^n(\Om)$, the
existence of a generalized solution, introduced in [Mi], was
investigated in [Gia, G2]. Here we consider solutions in
$SH_1(\Om)$. We say $u\in SH_1(\Om)$ is a {\it weak solution} of
(1.3) if $\mu_1[u]=\nu$.

Assume that for any Caccioppoli set $\omega\subset\Om$ with nonzero
measure,
$$\nu(\omega)<|\p \omega|. \tag 7.2$$
This is also a necessary condition for the existence of smooth
solutions to the mean curvature equation (7.1), which can be
verified easily by taking integration by parts of the equation.

Let $\rho$ be a mollifier, as was given in (5.3). Let $g_\eps(x)$
be the mollification of $\nu$, namely
$$g_\eps(x)=\int_{\Om}\rho_\eps(x-y)d\nu.$$
Extend $\nu$ to $\R^n$ such that $\nu=0$ outside $\Om$. Then
$g_\eps\in C^\infty(\R^n)$ and $g_\eps dx$ converges to $\nu$
weakly.

\proclaim{Lemma 7.1} For any open set $\omega\subset\Om$, we have
$$\int_{\omega}g_\eps dx<|\p\omega|. \tag 7.3$$
\endproclaim

\noo{\bf Proof}. We have
$$\align
\int_{\omega}g_\eps dx
 & = \int_{\omega}dx \int_{\Om}\rho_\eps (x-y) d\nu\\
 & = \int_{|z|\leq 1}\nu(\omega-\eps z)\rho(z)dz,\\
\endalign $$
where $\omega-\eps z=\{x\in\R^n:\ x+\eps x\in\omega\}$. By (7.2),
$\nu(\omega-\eps z)<|\p\omega|$. Hence we obtain (7.3). $\square$

Consider the approximation problem
$$ \align
H_1[u]& =g_\eps (x)\ \ \ \text{in}\ \ \Om,\tag 7.4\\
u & =\varphi\ \ \ \text{on}\ \ \pom.\\
\endalign$$
For equation with smooth right hand side, we quote the following
result [Gia].

\proclaim{Lemma 7.2} Under condition (7.3), there is a minimizer
$u_\eps$ of the functional
$$\Cal F(u)=\int_\Om \sqrt {1+|Du|^2}
  -\int_\Om g_\eps u+\int_{\pom} |u-\phi| .\tag 7.5$$
If $\phi\in C^0(\pom)$, the minimizer is a smooth solution to the
mean curvature equation (7.4). If the mean curvature $H'$ of
$\pom$ (with respect to the inner normal) satisfies
$$H'(x)>\frac{n}{n-1}g_\eps (x)\ \ \ \ \forall \ \ x\in\pom, \tag 7.6$$
then $u_\eps=\phi$ on $\pom$. \endproclaim

\noo{\bf Remark 7.1}. By our Harnack inequality, the minimizer is
a smooth solution to the mean curvature equation $H_1[u]=g_\eps$
if $\phi\in L^1(\pom)$.

\proclaim{Theorem 7.1} Let $\Om$ be a bounded domain in $\R^n$
with $C^2$ boundary. Let $\nu$ be a nonnegative measure which
satisfies (7.2) and can be decomposed as $\nu=\nu_1+f$ for some
nonnegative measure $\nu_1$ with compact support in $\Om$ and some
Lipschtiz function $f\ge 0$. Suppose the boundary mean curvature
satisfies
$$H'(x)>\frac{n}{n-1} f (x)\ \ \ \ \forall \ \ x\in\pom. \tag 7.7$$
Then (7.1) has a weak solution.
\endproclaim

\noo{\bf Proof}. We divide the proof into two steps.

\noo {\bf Step 1}. First we prove the theorem under the additional
assumption that there exists a positive constant $\eta>0$ such
that for any Caccioppoli set $\omega \subset\Om$,
$$\nu(\omega)\leq(1-\eta)|\p  \omega|.\tag 7.8$$
Let $g_\eps$ be the mollification of $\nu$ as above. Note that
(7.7) implies (7.6) for small $\eps>0$. Hence by Lemma 7.2, there
is a solution $u_\eps$ to (7.4). By Theorem 3.2, $u_\eps$ is
uniformly bounded,
$$\sup_{\pom}\phi\ge u_\eps\ge -C\tag 7.9$$
for some $C>0$ independent of $\eps$. By assumption, $\nu$ is
given by a Lipschitz continuous function $f$ in $\Om-S$, where
$S=\text{supp}\, \nu_1$. Hence $u_\eps$ is locally uniformly
bounded in $C^2(\Om-S)$. By Theorem 3.1, $u_\eps$ is uniformly
bounded in $W^{1,1}(\Om')$ for any $\Om'\Subset\Om$. Hence we may
assume that $u_\eps$ converges in $L^1$ to a limit function $u$.
Note that $g_\eps dx$ converges weakly to $\nu$. By Theorem 6.1,
$u$ is a weak solution of (7.1). By Corollary 3.1 and since
$\nu=f$ is Lipschitz continuous in $\Om-S$, $u\in BV(\Om)$.

\noo {\bf Step 2}. Next we remove the assumption (7.8). For any
small constant $\delta\in (0, 1)$, from Step 1 there is a solution
$u_\delta\in BV(\Om)$ to
$$ \align
H_1[u]& =(1-\delta)\nu \ \ \ \text{in}\ \ \Om,\tag 7.10\\
u & =\varphi\ \ \ \text{on}\ \ \pom.\\
\endalign$$
Then $u_\delta$ is monotone, namely $u_{\delta_1}\ge u_{\delta_2}$
if $\delta_1>\delta_2$; and $u_\delta$ is smooth near $\pom$. We
wish to prove that $u_\delta$ converges to a solution of (7.1) as
$\delta\to 0$. Since $u_\delta$ is monotone, we may define
$$u=\lim_{\delta\to 0} u_\delta. \tag 7.11$$
Denote $N=:\{x\in\Om:\ \ u(x)=-\infty\}$. If $N$ has measure zero,
then by Lemma 4.3, $u\in SH_1(\Om)$, and by Theorem 6.1,
$\mu_1[u]=\nu$. To see that $u$ satisfies the boundary condition
$u=\phi$ on $\pom$, note that $\nu=f$ is a Lipschitz function near
$\pom$ and recall that Lemma 2.5 holds for functions satisfying
$H_1[u]\ge f$, see Remark 2.4. Hence $u_\delta$ is locally uniformly
bounded and smooth near $\pom$. Hence the boundary condition
$u=\phi$ is satisfied and so $u$ is a weak solution of (7.1).

It remains to prove that Lebesgue measure $|N|=0$. Suppose to the
contrary that
$$|N|>\sigma>0.\tag 7.12$$
We claim that there exists a positive constant $\eta>0$ such that
$$\nu(\Om_t)<(1-\eta)|\pom_t|\tag 7.13$$
for all large $t$, where $\Om_t=\{x\in\Om:\ \ u(x)\le -t\}$, so
that $N=\Om_\infty$. (7.13) can be proved by a compactness
argument. Indeed, if it is not true, there is a sequence of
$\{t_j\}$, $t_j\to t_\infty\le \infty$, such that
$$\nu(\Om_{t_j})\ge (1-2^{-j})|\pom_{t_j}|.$$
Let $\phi_j$ be the characteristic function so that
$$|\pom_{t_j}|=\int_{\R^n}|D\phi_j|. $$
Since $\nu(\Om_{t_j})\le \nu(\Om)$ is uniformly bounded, $\phi_j$
converges in $L^1$ to the characteristic function $\phi$ of
$\Om_{t_\infty}$ and
$$\int_{\R^n}|D\phi|\le \lim_{j\to\infty}\int_{\R^n}|D\phi_j|.$$
Since $\Om_t$ is monotone, we have $\nu(\Om_{t_j})\to
\nu(\Om_{t_\infty})>\sigma$. Hence we obtain
$$\nu(\Om_{t_\infty})\ge \int_{\R^n}|D\phi|, $$
which is in contradiction with (7.2). Hence (7.13) holds.

Denote $\Om_{\delta, t}=\{x\in\Om:\ \ u_\delta(x)\le -t\}$. Recall
that $u_\delta$ is monotone. Hence for any $t>0$, $|\Om_{\delta,
t}|>\sigma$ provided $\delta$ is sufficiently small. For any fixed
$t$, by a compactness argument as above, we also have
$$\nu(\Om_{\delta, t})<(1- \eta)|\pom_{\delta, t}| \tag 7.14$$
when $\delta$ is sufficiently small.  Let $\delta_t>0$ be the sup
of all such $\delta$. Then again by a similar compactness
argument, we have
$${\underline\lim}_{t\to t_0}
                    \delta_{t}\ge \delta_{t_0}. \tag 7.15$$
Therefore for any $T>0$, we can choose $\delta>0$ sufficiently
small such that (7.14) holds for all $t\in (0, T]$. Now we fix $T$
as in (3.5).  By Step 1 above, $u_\delta$ is a bounded function
and $u_\delta\in BV(\Om)$. Hence the proof of Theorem 3.2 is valid
(see Remark 3.1) and we obtain
$$\inf u_\delta\ge -C$$
for some $C>0$ depending on $n$, $|\Om|$, $\inf_{\pom} u_\delta$,
and $\eta$, but is independent of $\delta$. Sending $\delta\to 0$,
we find that $u$ is bounded from below, a contradiction. $\square$

\vskip10pt

\noo{\bf Remark 7.2}. Condition like (7.2) was included in [Gia,
G1, G2]. When $\nu$ (more precisely its density) is a bounded
function, (7.2) implies (7.8) for a small $\eta$ [G1].

\noo{\bf Remark 7.3}. A weak solution is usually not $C^2$ smooth
if $\nu$ is not Lipschitz continuous. This is easily seen by
considering functions of one variable, $u=u(x_1)$. However, if
$n\le 7$ and $\nu$ is a bounded function and the weak solution is
a minimizer of the functional (7.5),  then the graph of the
solution is a $C^{2, \alpha}$ hypersurface if $\nu$ is a H\"older
continuous function; or $C^{1, \alpha}$ if $\nu$ is a bounded
nonnegative function [Ma].

\vskip20pt

\centerline{\bf 8. Remarks}

 \vskip10pt

We include an example showing that some potential theoretical
properties which hold for the $p$-Laplace equation and the
$k$-Hessian equation [HKM, Lab, TW1-TW4] may not hold for
curvature equations.

Let
$$ u_c(x)= \cases
   a(r-1)^\delta\ \ &\text{if}\ \ r\ge 1,\\
   -b(1-r)^\sigma-c\ \ &\text{if}\ \ 0\le r<1,\\
   \endcases$$
where $r=|x|$, $a, b>1, c\ge 0, \delta, \sigma\in (0, \frac 12)$
are positive constants.  Then $H_1[u_0]$ is positive and H\"older
continuous near $\p B_1$,  but $u_0\not\in C^1$, since
$|Du|=\infty$ on the sphere $\{|x|=1\}$. As remarked at the end of
last section, the graph of $u_0$ is $C^{2, \alpha}$ for some
$\alpha>0$.

If $c>0$, $u_c$ is $H_1$-subharmonic, and can be approximated by
smooth $H_1$-subharmonic functions. Therefore weak solutions to
the Dirichlet problem (7.1), without the restriction (7.2), is not
unique in general. We note that the corresponding uniqueness
problem for the $p$-Laplace equation and the $k$-Hessian equations
remains open.

When $c>0$, we also see that the Wolff potential estimate (see,
e.g., [L, TW4]) does not hold for the mean curvature equation, and
an $H_1$-subharmonic function may not be quasi-continuous, as the
capacity of $\p B_1$ is positive.

 \vskip20pt

\baselineskip=12.0pt
\parskip=2.0pt

\Refs\widestnumber\key{ABC}

\item {[AS]}  L. Ambrosio and H.M. Soner,
       Level set approach to mean curvature flow in arbitrary codimension,
       J. Diff. Geom. 43 (1996), 693--737.

\item {[FL]} R. Finn and J. Lu,
       Regularity properties of H-graphs,
       Comment. Math. Helv., 73(1998), 379-399.

\item {[F1]} H. Federer, Curvature measures.
       Trans. Amer. Math. Soc. 93(1959), 418--491.

\item {[F2]} H. Federer,
       Geometric measure theory,
       Springer,  New York, 1969.

\item {[GT]} D. Gilbarg and N.S. Trudinger,
       Elliptic Partial Differential Equations of Second Oder,
       Springer, Second Edition, 1983.

\item {[Gia]} M. Giaquinta,
       On the Dirichlet problem for surfaces of prescribed mean curvature,
       Manuscripta Math. 12 (1974), 73--86.

\item {[G1]} E. Giusti,
       On the equation of surfaces of prescribed mean curvature,
       Existence and uniqueness without boundary conditions,
       Invent. Math. 46(1978), 111-137.

\item {[G2]} E. Giusti,
       Generalized solutions for the mean curvature equation,
       Pacific J. Math., 88(1980), 297-321.

\item {[G3]} E. Giusti,
       Minimal surfaces and functions of bounded variations,
       Birkhauser, Boston, 1984.

\item {[HKM]} J. Heinonen, T. Kilpelainen and O. Martio,
       Nonlinear Potential Theory of Degenerate Elliptic Equations,
       Oxford Univ. Press, 1993.

\item {[K1]} N. Korevaar,
       An easy proof of the interior gradient bound for solutions
       to the prescribed mean curvature equation,
       In {\it Nonlinear functional analysis and its applications},
       Proc. Sympos. Pure Math., 45 (Part 2), 81--89,
       Amer. Math. Soc., 1986.

\item {[K2]} N. Korevaar,
       A priori interior gradient bounds for solutions to elliptic Weingarten
       equations, Ann. Inst. H. Poincar\'e Anal. Non Lin\'eaire, 4(1987),
       405--421.

\item {[Lab]} D.A. Labutin,
       Potential estimates for a class of fully nonlinear elliptic
       equations, Duke Math. J. 111 (2002), 1--49.

\item {[Lia]}  F.T. Liang,
       Harnack's inequalities for solutions to the mean curvature
       equation and to the capillarity problem,
       Rend. Sem. Mat. Univ. Padova., 110(2003), 57-96.

\item {[Ma]} U. Massari,
       Esistenza e regolarita delle ipersuperfice di curvatura media
       assegnata in $R\sp{n}$ (Italian),
       Arch. Rat. Mech. Anal. 55 (1974), 357--382.

\item {[Mi]} M. Miranda,
       Superficie minime illimitate (Italian),
       Ann. Scuola Norm. Sup. Pisa Cl. Sci. (4) 4(1977), 313--322.

\item {[PS1]} P. Pucci and J. Serrin,
       The Harnack inequality in $R^2$ for quasilinear elliptic
       equations, J. Anal. Math., 85(2001), 307-321.

\item {[PS2]} P. Pucci and J. Serrin,
       The strong maximum principle revisited,
       J. Diff. Eqns 196 (2004),  1--66.

\item {[San]} L.A. Santal\'o,
       Integral geometry and geometric probability,
       Addison-Wesley, 1976.

\item {[S]} R. Schneider,
       Convex bodies: the Brunn-Minkowski theory,
       Cambridge Univ. Press, Cambridge, 1993.

\item {[T1]} N.S. Trudinger,
        Harnack inequalities for nonuniformly elliptic divergence
        structure equations,
        Invent. Math. 64 (1981), 517--531.

\item {[T2]} N.S. Trudinger,
       A priori bounds and necessary conditions for solvability of
       prescribed curvature equations,
       Manuscripta Math. 67(1990), 99-112.

\item {[TW1]} N.S. Trudinger and X.-J. Wang,
       Hessian measures I,
       Topol. Methods Nonlinear Anal., 10(1997), 225-239.

\item {[TW2]} N.S. Trudinger and X.-J. Wang,
       Hessian measures II, Ann. Math., 150(1999), 579-604.

\item {[TW3]} N.S. Trudinger and X.-J. Wang,
        Hessian measures III,
        J. Funct. Anal., 193 (2002), 1-23.

\item {[TW4]} N.S. Trudinger and X.-J. Wang,
       On the weak continuity of elliptic operators and applications
       to potential theory, Amer. J. Math., 124(2002), 369-410.

\item {[Wan]} X.-J. Wang,
       Interior gradient estimates for mean curvature equations,
       Math. Z. 228 (1998),  73--81.

\endRefs

\enddocument

\end